\documentclass[12pt]{article}
\usepackage{latexsym}
\usepackage{amssymb,amsmath}
\usepackage{epsfig}
\usepackage{amsfonts}

\newtheorem{theorem}{Theorem}[section]

\newtheorem{corollary}[theorem]{Corollary}
\newtheorem{proposition}[theorem]{Proposition}

\newcommand{\proof}{\noindent{\bf Proof.\ }}
\newcommand{\qed}{\hfill $\square$ \bigskip}
\def\cp{\,\square\,}
\newcommand{\dl}{{\rm dim_{\ell}}}

\begin{document}

\title{Local metric dimension of graphs: generalized hierarchical products and some applications}

\author{Sandi Klav\v zar $^{a,}$\footnote{The financial support of the Slovenian Research Agency (research core funding No.\ P1-0297 and project J1-9109) is acknowledged.} 
\and Mostafa Tavakoli $^{b}$} 

\date{}

\maketitle

\begin{center}
$^a$ Faculty of Mathematics and Physics, University of Ljubljana, Slovenia\\
{\tt sandi.klavzar@fmf.uni-lj.si}

\medskip

$^d$ Department of Applied Mathematics, Faculty of Mathematical Sciences,\\
Ferdowsi University of Mashhad, P.O.\ Box 1159, Mashhad 91775, Iran\\
{\tt m$\_$tavakoli@um.ac.ir}

\end{center}

\begin{abstract}
Let $G$ be a graph and $S\subseteq V(G)$. If every two adjacent vertices of $G$ have different metric $S$-representations, then $S$ is a local metric generator for $G$. A local metric generator of smallest order is a local metric basis for $G$, its order is the local metric dimension of $G$. Lower and upper bounds on the local metric dimension of the generalized hierarchical product are proved and demonstrated to be sharp. The results are applied to determine or bound the dimension of several graphs of importance in mathematical chemistry. Using the dimension, a new model for assigning codes to customers in delivery services is proposed. 
\end{abstract}


\medskip\noindent
{\bf AMS Subj.\ Class:} 05C12; 05C76; 90B80

\section{Introduction}
\label{sec:intro}

All graphs considered in this paper are connected and simple. If $G= (V(G), E(G)$ is a graph, then its order and its size are denoted with $n(G)$ and $m(G)$, respectively. If $u,v\in V(G)$, then $d_G(u,v)$ denotes the standard shortest-path distance between $u$ and $v$ in $G$, that is, the number of edges on a shortest $u,v$-path. If $S=\{v_1,\ldots, v_k\}\subseteq V(G)$, then the {\em metric $S$-representation} of a vertex $v\in V(G)$ is the vector 
$$r_G(v|S) = (d_G(v,v_1), \ldots, d_G(v,v_k))\,.$$ 
A {\em metric generator} for $G$ is a vertex subset $S\subseteq V(G)$ such that the vertices of $G$ have pairwise different metric $S$-representations. A metric generator of smallest order is a {\em metric basis} for $G$, its order being the {\em metric dimension} ${\rm dim}(G)$ of $G$. 

The metric dimension was introduced in~\cite{Harary, Slater} and is used to model many real world problems. These include navigation of robots~\cite{khuller-1996} and chemical problems~\cite{johnson-1993}. However, quite often we do not need to distinguish all pairs of vertices but only adjacent ones. From this reason the local metric dimension was introduced in~\cite{okamoto-2010}. To be more specific, $S\subseteq V(G)$ is a {\em local metric generator} for $G$ if the condition of having different metric $S$-representations is fulfilled for every {\bf adjacent} vertices of $G$. The {\em local metric dimension} ({\em lmd} for short) $\dl(G)$ of $G$ is then, as expected, the smallest order of a local metric generator, and such a set is a {\em local metric basis} ({\em lmb} for short) for $G$. We mention here two further, recently proposed variants of the metric dimension. In~\cite{hakanen-2018} resolving sets locate up to some fixed $\ell$, $\ell \ge 1$, vertices in a graph, while in~\cite{kelenc-2018} resolving sets locate the edges of a graph. The property of being located is often also combined with some other properties, say being a dominating set; for a comparison of resolving sets with locating-dominating set and additional related sets see~\cite{gonzales-2018}. 

The lmd has been by now investigated on different graph operations. Already in the seminal paper~\cite{okamoto-2010} it was considered on the Cartesian product. Rodr\'{i}guez-Vel\'{a}zquez, Garc\'{i}a G\'{o}mez, and Barrag\'{a}n-Ram\'{i}rez followed with studies of it on rooted products~\cite{rooted}. The behavior of the lmd on corona products and on edge corona produces was investigated by Rodr\'{i}guez-Vel\'{a}zquez, Barrag\'{a}n-Ram\'{i}rez, and Garc\'{i}a G\'{o}mez~\cite{Corona} and by Rinurwati, Slamin, and Suprajitno~\cite{ecorona}, respectively. Finally Barrag\'an-Ram\'irez and Rodr\'iguez-Vel\'azquez studied it under the strong product operation~\cite{strong}, while in~\cite{fernau-2018} Fernau and Rodr\'{i}guez-Vel\'{a}zquez related the local dimensions of corona and strong products with the newly introduced adjacency metric dimension. We also refer to~\cite{saputro-2016} for the lmd of regular graphs. 

The aim of the present work is to continue the above line of investigation on graph operations by considering the lmd of the generalized hierarchical product (to be defined below). We refer to~\cite{anderson-2018, hossein-2017} as well as references therein for results on this graph operation. The generalized hierarchical product is in particular important because it generalizes several classical operations, such as the Cartesian product. 

In the next section we give lower and upper bounds on the lmd of the generalized hierarchical product and demonstrate their sharpness. Some earlier results are shown to be consequences of the present results. Then, in Section~\ref{sec:applications}, we first apply the results of Section~\ref{sec:hierarchical} to determine or bound the lmd of several graphs that are important in mathematical chemistry. We conclude the paper with a new model for assigning codes to customers in delivery services. The model uses local metric bases and is in many cases significantly more efficient that the classical model from~\cite{khuller-1996}.   

In the rest of the section we formally introduce the generalized hierarchical product and recall the key result about the distances in it. Let $G$ and $H$ be graphs and let $\emptyset \ne U\subseteq V(G)$. The {\em generalized hierarchical product} $G(U)\sqcap H$ (with respect to $U$) of $G$ and $H$, is a graph with the vertex set  
$$V(G(U)\sqcap H) = V(G)\times V(H)\,,$$ 
and the edge set  
$$\{(g,h)(g',h'):\ gg'\in E(G), h=h'\} \cup \{(g,h)(g',h'):\ g=g'\in U, hh'\in E(H)\}\,,$$
see~\cite{barriere-2009}. 

If $U\subseteq V(G)$ and $u,v\in V(G)$, then we say that a $u,v$-walk $W$ is a {\em $u,v$-walk through $U$} if $W$ is an $u,v$-walk in $G$ that contains some vertex of $U$. (Note that this vertex from $U$ could be one of $u$ and $v$.) With $d_{G(U)}(u,v)$ we denote the length of a shortest $u,v$-walk through $U$. The following fundamental observation from~\cite{barriere-2009} will be used throughout the paper, mostly without explicitly mentioning it.

\begin{proposition} 
\label{prp:distance}
If $G$ is a graph with $U\subseteq V(G)$ and $H$ is a graph, then 
$$d_{G(U)\sqcap H}((g,h),(g',h'))=\begin{cases}
d_{G(U)}(g,g')+d_H(h,h'); &  h\neq h',\\
d_G(g,g'); & h=h'.
\end{cases}$$
\end{proposition}

\section{Generalized hierarchical products}
\label{sec:hierarchical}

It this section we first prove bounds on the studied dimension of generalized hierarchical product $G(U)\sqcap H$ under the condition that $G$ contains a lmb which is contained in $U$. Although this is not always the case, the result has several interesting consequences. In particular, if $G$ is bipartite, then $\dl(G\cp H) = \dl(H)$. In our second main result we give a general upper bound on the lmd of $G(U)\sqcap H$ and show that it is also sharp. To formulate the first main result we need to extend the concept of the lmd as follows. 

Let $G$ be a graph and let $\emptyset \ne U\subseteq V(G)$. If $S=\{v_1,\ldots, v_k\}\subseteq V(G)$, then we say that an {\em $U$-metric $S$-representation} of $v\in V(G)$ is  
$$r_{G(U)}(v|S) = (d_{G(U)}(v,v_1), \ldots, d_{G(U)}(v,v_k))\,.$$ 
$S\subseteq V(G)$ is a {\em $U$-metric local generator} if every two end-points of an edge of $G$ have pairwise different $U$-metric local representations. A $U$-metric local generator of smallest order is a {\em $U$-metric local basis} for $G$. The number of vertices in it is the {\em $U$-metric local dimension} $\dl(G|U)$ of $G$. In the case when $U = V(G)$, we have $\dl(G|U)$ = $\dl(G)$. 

\begin{theorem}
\label{thm:basis-under-U}
Let $G$ and $H$ be graphs and $\emptyset \ne U\subseteq V(G)$. If $G$ contains a lmb which is contained in $U$, then 
$$ \max\{\dl(G|U), \dl(H)\} \le \dl(G(U)\sqcap H)\leq \max\{\dl(G), \dl(H)\}\,.$$
\end{theorem}

\proof
Set $X = G(U)\sqcap H$ for the rest of the proof. Let $S_G = \{g_1, \ldots, g_{\dl(G)}\}$ be a lmb of $G$ that is contained in $U$, and let $S_H = \{h_1, \ldots, h_{\dl(H)}\}$ be a lmb of $H$. 

We first show that $\dl(X) \le \max\{\dl(G), \dl(H)\}$ and for this sake consider the following cases. 

\medskip\noindent
{\bf Case 1}: $k = \dl(G) \ge \dl(H)$. \\
Set $S=\{(g_i,h_i):\ i\in [k]\}$, where the indices in the second coordinate are taken modulo $\dl(H)$. We claim that $S$ is a local metric generator for $X$ and for this sake consider arbitrary adjacent vertices $(g, h)$ and $(g',h')$ in $V(X)\setminus S$.

Suppose first that $h=h'$. Then in $G$ there exists $g_i\in S_G$ such that $d_G(g,g_i)\ne d_G(g',g_i)$. Since $g_i\in U$ it follows that $d_{G[U]}(g,g_i) = d_G(g,g_i) \ne d_G(g',g_i) = d_{G[U]}(g',g_i)$. (It is possible that $g_i = g$ or $g_i = g'$.) Applying Proposition~\ref{prp:distance} we thus infer that
\begin{eqnarray*}
d_{X}((g,h), (g_i,h_i)) & = & d_{G(U)}(g,g_i) + d_H(h,h_i) \\
& \ne & d_{G(U)}(g',g_i) + d_H(h,h_i) \\
& = & d_{X}((g',h), (g_i,h_i))\,.
\end{eqnarray*}
It follows that $r_X((g,h)|S) \ne r_X((g',h')|S)$. 

Suppose second that $g=g'$. Since $hh'\in E(H)$ we have $h_i\in S_H$ with $d_H(h,h_i)\ne d_H(h',h_i)$. Then e  
\begin{eqnarray*}
d_{X}((g,h), (g_i,h_i)) & = & d_{G(U)}(g,g_i) + d_H(h,h_i) \\
& \ne & d_{G(U)}(g,g_i) + d_H(h',h_i) \\
& = & d_{X}((g,h'), (g_i,h_i))\,,
\end{eqnarray*}
which in turn implies that $r_X((g,h)|S) \ne r_X((g',h')|S)$. 

It follows from the above that $\dl(X) \le |S| = k = \max\{\dl(G), \dl(H)\}$.

\medskip\noindent
{\bf Case 2}: $k = \dl(H) \ge \dl(G)$. \\
In this case set again $S=\{(g_i,h_i):\ i\in [k]\}$, except that now the indices in the first coordinate are taken modulo $\dl(G)$. The proof now proceeds analogously as in Case~1. Therefore, also in this case we have $\dl(X) \le |S| = k = \max\{\dl(G), \dl(H)\}$. From Cases~1 and 2 we conclude that $\dl(X) \le \max\{\dl(G), \dl(H)\}$. 

\medskip
Let now $S$ be an arbitrary lmb for $X$. Consider first the projection $S_H$ of $S$ on $H$, that is, $S_H = \{h\in V(H):\ \exists g\in V(G)\ {\rm such\ that}\ (g,h)\in S\}$. Let $h, h'\in V(H)\setminus S_H$ such that $hh' \in E(H)$. Let $g\in U$. Then $(g,h)(g,h')\in E(X)$ and $(g,h), (g,h')\in V(X)\setminus S$. As $S$ is a lmb for $X$, there is a vertex $(g'',h'')\in S$ with $d_X((g,h), (g'',h''))\ne d_X((g,h'), (g'',h''))$. Since 
$$d_{X}((g,h), (g'',h'')) = d_{G(U)}(g,g'') + d_H(h,h'')$$ 
and 
$$d_{X}((g,h'), (g'',h'')) = d_{G(U)}(g,g'') + d_H(h',h'')\,,$$ 
it follows that $d_H(h,h'') \ne d_H(h',h'')$. As $h''\in S_H$ it follows that $S_H$ is a local metric generator of $H$ which means that $\dl(H) \le |S_H| \le |S| = \dl(X)$. 

Consider second the projection $S_G$ of $S$ on $G$, that is, $S_G = \{g\in V(G):\ \exists h\in V(H)\ {\rm such\ that}\ (g,h)\in S\}$. Let $g, g'\in V(G)\setminus S_G$ such that $gg' \in E(G)$. Let $h\in V(H)$. Then $(g,h)(g',h)\in E(X)$ and $(g,h), (g',h)\in V(X)\setminus S$. As $S$ is a lmb for $X$, there exists a vertex $(g'',h'')\in S$ such that $d_X((g,h), (g'',h''))\ne d_X((g',h), (g'',h''))$. Since 
$$d_{X}((g,h), (g'',h'')) = d_{G(U)}(g,g'') + d_H(h,h'')$$ 
and 
$$d_{X}((g',h), (g'',h'')) = d_{G(U)}(g',g'') + d_H(h,h'')\,,$$ 
it follows that $d_{G(U)}(g,g'') \ne d_{G(U)}(g',g'')$. As $h''\in S_G$ it follows that $S_G$ is a $U$-metric local generator of $G$ and so $\dl(G|U) \le |S_G| \le |S| = \dl(X)$.
\qed

Theorem~\ref{thm:basis-under-U} implies several exact results. For instance, if $U = V(G)$ then $G(U)\sqcap H$ is just the Cartesian product $G\cp H$ and $\dl(G|U) = \dl(G)$, hence we get:

\begin{corollary} {\rm \cite{okamoto-2010}}
\label{cor:Cartesian}
For any $G$ and $H$, $\dl(G\cp H) = \max \{\dl(G), \dl(H)\}$.
\end{corollary}

Note that if $G$ is bipartite, then $\dl(G|U) = 1$ holds for any $\emptyset \ne U\subseteq V(G)$. Therefore: 

\begin{corollary}
\label{cor:bipartite}
If $G$ is a bipartite graph and $H$ a graph, then $\dl(G\cp H) = \dl(H)$.
\end{corollary}

The {\em join} $G+H$ of disjoint graphs $G$ and $H$ is obtained from their disjoint union by adding all possible edges between the vertices from $G$ and the vertices from $H$. The {\em corona} $G\odot H$ of graphs $G$ and $H$ is obtained from the disjoint union of a copy of $G$ and $n(G)$ copies of $H$, where each vertex of the $i^{\rm th}$ copy of $H$ is adjacent to the $i^{\rm th}$ vertex of $G$, $i\in [n(G)]$. Note that $G\odot H = (H+K_1)(U)\sqcap G$ where $U = V(K_1)$. Among many results from~\cite{Corona} we extract~\cite[Corollary 5(i)]{Corona} which asserts that 
$$\dl(H\oplus K_t) = n(H)(t-1)\,,$$
that is,  
$$\dl(K_{t+1}(\{v\}) \sqcap H)) = n(H)(t-1)\,.$$
This result demonstrates that the assumption of Theorem~\ref{thm:basis-under-U} that the first factor must contain a lmb which is contained in $U$ cannot be avoided. Since in general this condition is not fulfilled, we state the following bound for the general case. 

\begin{theorem}
\label{thm:general-case}
Let $G$ and $H$ be graphs and $\emptyset \ne U\subseteq V(G)$. If $S_G$ is a lmb of $G$ such that $|S_G\cap U| = k$, then
$$\dl(G(U)\sqcap H)\leq n(H)(\dl(G) - k) + k\,.$$
\end{theorem}

\proof
Set $X = G(U)\sqcap H$ for the rest of the proof. Let $S_G=\{g_1, \ldots, g_{\dl(G)}\}$ and let $S_H=\{h_1,\ldots, h_{\dl(H)}\}$ be a lmb of $H$. Assume w.l.o.g.\ that $S_G\cap U=\{g_1, \ldots, g_k\}$. We claim that 
$$S=\big((S_G\setminus U)\times V_H\big)\cup \{(g_i,h_i):\ i\in [k]\}\,,$$
is a local metric generator for $X$, where the indices $i$, if necessary, are taken modulo $k$. For this sake consider a pair of adjacent vertices $(g, h)$ and $(g',h')$ from $V(X)\setminus S$ and distinguish the following two natural cases further divided into subcases. 

\medskip\noindent
{\bf Case 1}: $gg'\in E(G)$ and $h=h'$. 

\medskip\noindent
{\bf Case 1.1}: $g, g'\notin S_G$. \\
Then there exists $g_i\in S_G$ such that $d_G(g_i,g)\neq d_G(g_i,g)$. If $i > k$, then $(g_i,h)\in S$ and it follows immediately that $d_{X}((g_i,h),(g,h))\ne d_{X}((g_i,h),(g',h))$. On the other hand, if $i\le k$, then since $g_i\in U$ we have $d_{G(U)}(g_i,g)\neq d_{G(U)}(g_i,g')$, which again implies that $d_{X}((g_i,h),(g,h))\ne d_{X}((g_i,h),(g',h))$.

\medskip\noindent
{\bf Case 1.2}: $g, g'\in S_G$. \\
In this subcase $g$ and $g'$ must both be from $S_G\cap U$. Hence there exists an $i$ such that $g_i=g$ and therefore, $|d_{X}((g_i,h_i),(g,h))-d_{X}((g_i,h_i),(g',h))|= d_G(g_i,g') = 1$. 

\medskip\noindent
{\bf Case 1.3}: $g\in S_G$, $g'\notin S_G$. \\
Since $g\in S_G$, we have $g=g_i$ for some $i\in [\dl(G)]$. If $i>k$, then 
$$d_{X}((g_i,h),(g,h))-d_{X}((g_i,h),(g',h))= -d_G(g_i,g')\neq 0,$$
and if $i\leq k$, then
$$d_{X}((g_i,h_i),(g,h))-d_{X}((g_i,h_i),(g',h))= -d_{G(U)}(g_i,g')\neq 0\,,$$
where again the index $i$ is taken modulo $k$ if necessary. 

\medskip\noindent
{\bf Case 2}: $g = g'\in U$ and $hh'\in E(H)$. 

\medskip\noindent
{\bf Case 2.1}: $g\notin S_G$ and $h,h'\notin S_H$. \\
In this case there exists $h_j\in S_H$ such that $d_H(h_j, h)\neq d_H(h_j,h')$ and so 
$$d_{X}((g_i,h_j),(g,h))-d_{X}((g_i,h_j),(g,h'))=d_H(h_j,h))-d_H(h_j,h')\neq 0\,,$$
where $i\in [\dl(G)]$.

\medskip\noindent
{\bf Case 2.2}: $h$ or $h'$ is in $S_H$ and $g\notin S_G$. \\
Without loss of generality we may assume that $h\in S_H$. Thus $d_{X}((g_i,h),(g,h))-d_{X}((g_i,h),(g,h'))=-d_H(h,h')\neq 0$, where $i>k$.

\medskip\noindent
{\bf Case 2.3}: $g\in S_G$ and $h,h'\notin S_H$. \\
In this subcase, $g$ must be in $S_G\cap U$. Thus, there exists $h_j\in S_H$ such that $d_H(h_j, h)\neq d_H(h_j,h')$ and so $d_{X}((g_i,h_j),(g,h))-d_{X}((g_i,h_j),(g,h'))=d_H(h_j,h))-d_H(h_j,h')\neq 0$, where again $i\in [\dl(G)]$.

In conclusion, for every pair of adjacent vertices $(g, h)$ and $(g',h')$ from $V(X)\setminus S$ there exists a vertex $(g_i, h_j)$ in $S$ such that $d_{X}((g_i,h_j),(g,h))\neq d_{X}((g_i,h_j),(g',h))$. Consequently, since clearly $|S| =  n(H)(\dl(G) - k) + k$ holds, the argument is complete. 
\qed

As already observed above, $G\odot H = (H+K_1)(U)\sqcap G$ where $U = V(K_1)$. If $H$ has radius more than $3$, then $v$ is not an element of any lmb of $H+K_1$. Consequently, by Theorem~\ref{thm:general-case}, we have 
$$\dl(G\odot H)=\dl((H+K_1)(U)\sqcap G)\leq n(G)\dl(H+K_1)\,.$$
On the other hand, it was proved in~\cite{Corona} that if $G$ is a connected graph and $H$ is a graph of radius at least $4$, then $\dl(G\odot H)=n(G)\dl(K_1+H)$. This demonstrates the sharpness of the bound of Theorem~\ref{thm:general-case}.

\section{Applications}
\label{sec:applications}

In this section we consider a couple of applications of the lmd. In the first part we apply the results of Section~\ref{sec:hierarchical} to determine or bound the lmd of some graphs that are important in mathematical chemistry. In the subsequent subsection we modify the application of the metric dimension from~\cite{khuller-1996} to delivery services such that the lmd is involved. In the new  model the length of constructed codes is lowered. 

\subsection{Local metric dimension of some molecular graphs}

Molecular graphs are, roughly speaking, graphs with the largest degree at most $4$. Several graph invariants of relevance in mathematical chemistry have already been investigated on the generalized hierarchical product, cf.~\cite{eliasi-2011, patt-2012}. On the other hand, the lmd has been used in~\cite{chem1, chem2} to classify certain molecules (molecular graphs). Here we add to these results additional examples of the metric dimension of molecular graphs. The latter can be represented as generalized hierarchical products, which in turn makes applicable the results of Section~\ref{sec:hierarchical}. 

\begin{figure}[t!]
\centerline{\includegraphics[scale=.35]{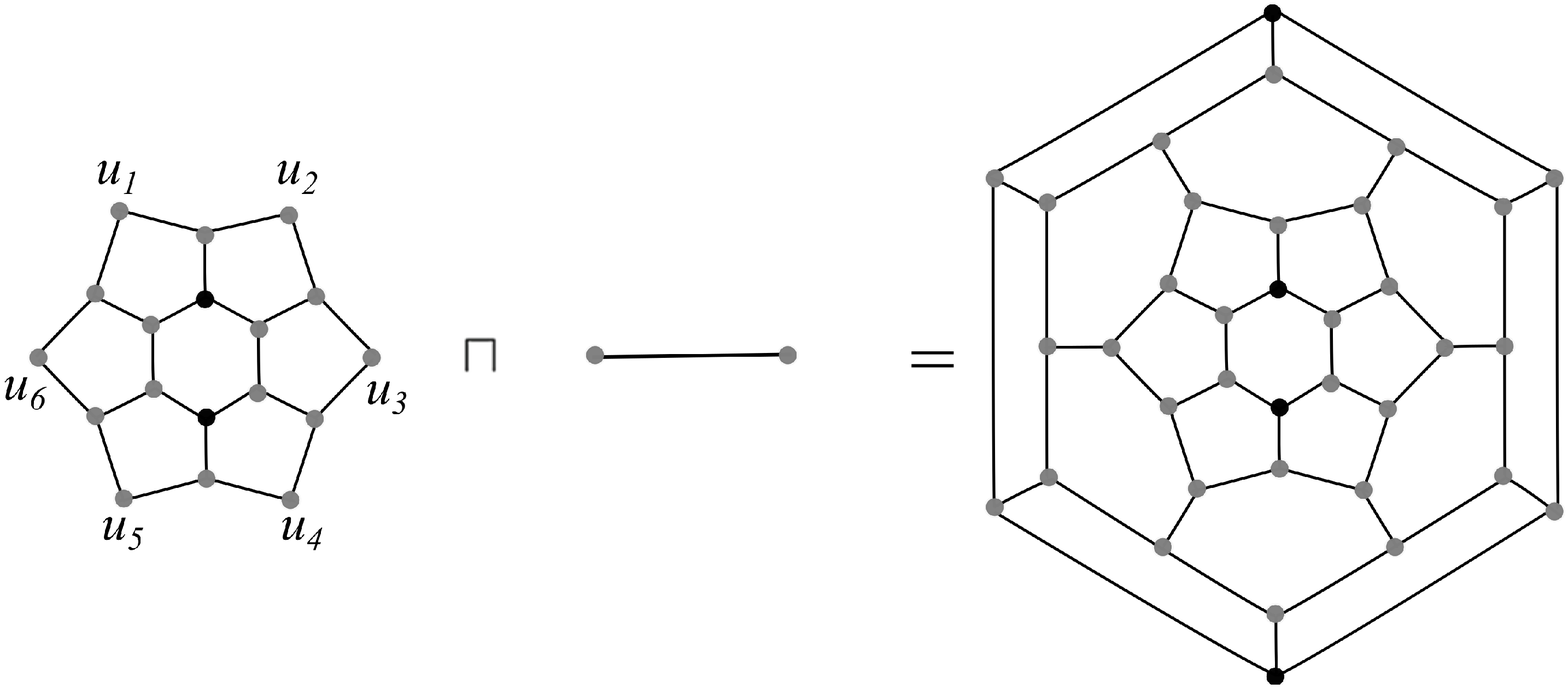}}
\caption{$G(U)\sqcap P_2=F_{5,12}$ where $U=\{u_1,\ldots, u_6\}$.}
\label{fig:F56}
\end{figure}

\medskip\noindent{\bf Example 1.} 
Consider the graph $G$ and its vertex subset $U=\{u_1,\ldots, u_6\}$ as shown on Fig.~\ref{fig:F56} (left). Then the generalized hierarchical product  $G(U)\sqcap P_2$ is the fullerene graph denoted $F_{5,12}$ drawn on Fig.~\ref{fig:F56} (right). 

As $G$ is not bipartite, $\dl(G)\ge 2$. On the other hand, one can check that the two black vertices of $G$ from the figure form a local metric generator. Consequently, $\dl(G)=2$. Then, using Theorem~\ref{thm:general-case}, $\dl(F_{5,12})=\dl(F_{5,6}(U)\sqcap P_2)\leq 4$. (The black vertices of $F_{5,12}$ from the figure form a local metric generator.)

\begin{figure}[t!]
\centerline{\includegraphics[scale=.27]{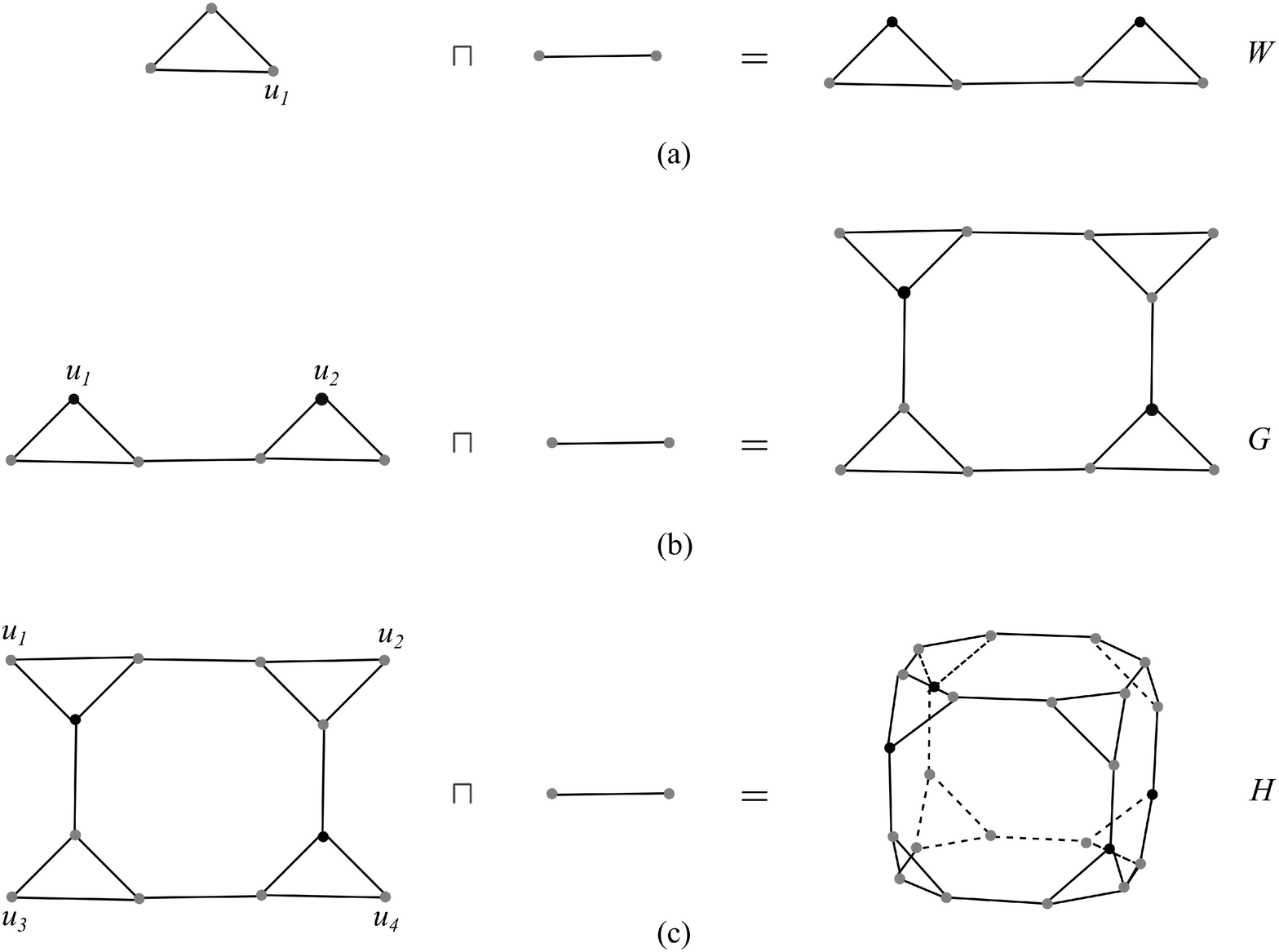}}
\caption{The molecular graph $H$ constructed via the generalized hierarchical product.}
\label{fig:truncated}
\end{figure}

\medskip\noindent
{\bf Example 2.} 
Let $H$ be the graph obtained from the of the truncated cube, see Fig.~\ref{fig:truncated} (bottom right). In the same figure it is shown how $H$ can be constructed in three steps from the triangle each time using the generalized hierarchical operation. Going backwards, $H=G(U)\sqcap P_2$, where $U=\{u_1, u_2, u_3, u_4\}$, see Fig.~\ref{fig:truncated}(c). Further, $G=W(U)\sqcap P_2$, where $\{u_1, u_2\}$, see Fig.~\ref{fig:truncated}(b). Finally (or firstly), $W=C_3(U)\sqcap P_2$, where $U=\{u_1\}$, see Fig.~\ref{fig:truncated}(a). 

The two black vertices of $W$ from the figure form a local metric generator of $W$. Consequently, $\dl(W)= 2$. In Fig.~\ref{fig:truncated}(b), $U$ is a local metric generator for $W$ and so by Theorem~\ref{thm:basis-under-U} we get $\dl(G)= \dl(W(U)\sqcap P_2)\leq 2$ and hence $\dl(G)=2$. Finally, none of the black vertices are in $U$, hence $k = 0$ in Theorem~\ref{thm:general-case} and so $\dl(H)= \dl(G(U)\sqcap P_2)\leq 4$.

\begin{figure}[t!]
\centerline{\includegraphics[scale=.3]{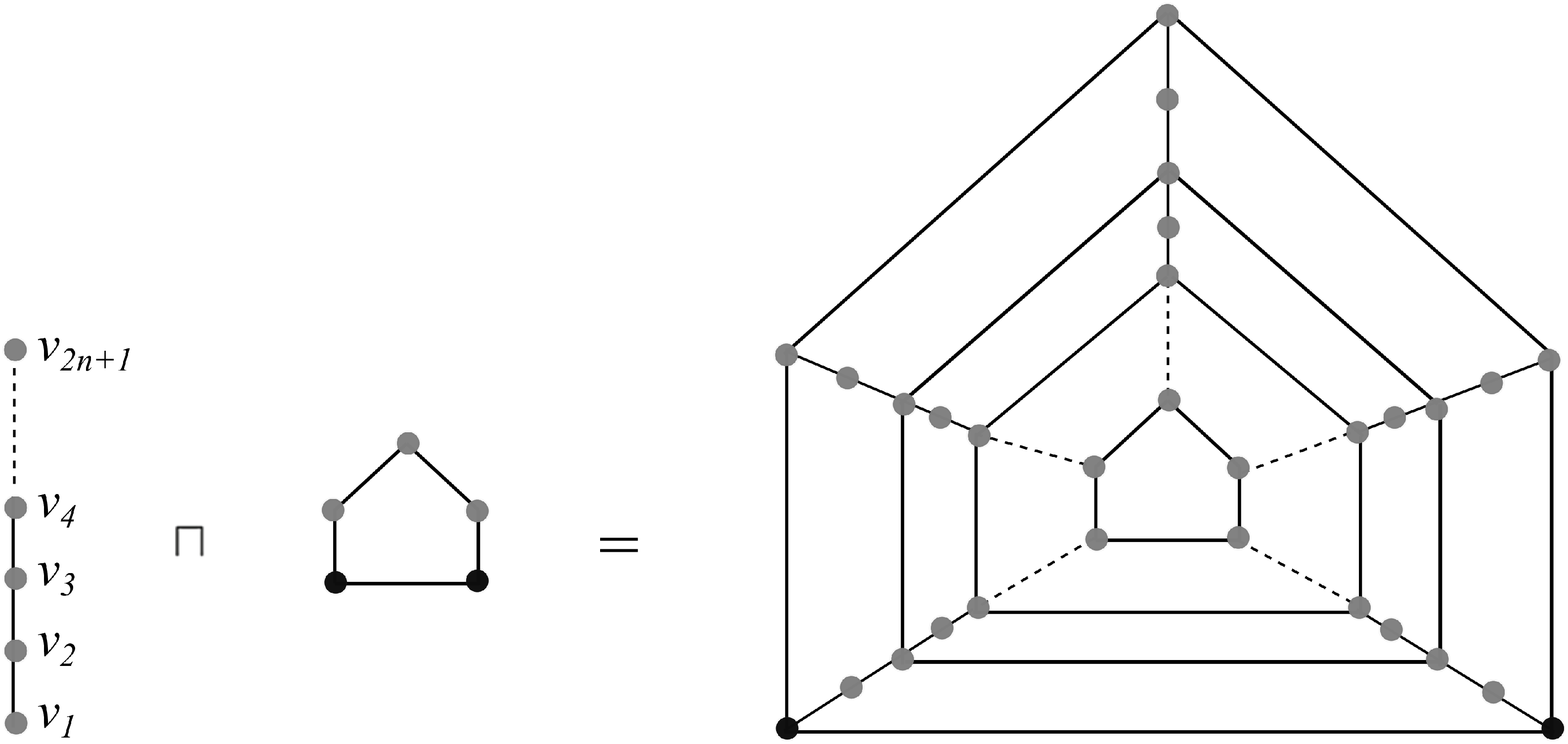}}
\caption{$\Gamma_{n,5}$ represented as a generalized hierarchical product.}
\label{fig:56-graph}
\end{figure}

\medskip\noindent
{\bf Example 3.} 
Let the vertices of the path $P_n$ be $v_1,\ldots, v_n$ in the natural order. Then consider the generalized hierarchical products $\Gamma_{n,k} = P_{2n+1}(U)\sqcap C_k$, where $k\ge 3$ and $U=\{v_{2i+1}:\ 0\leq i \leq k\}$. Fig.~\ref{fig:56-graph} displays the construction of $\Gamma_{n,5}$. 

If $k$ is an even number, then $\Gamma_{n,k}$ is bipartite and therefore $\dl(\Gamma_{n,k})=1$. For the other case we have: 

\begin{proposition}
If $k\ge 3$ is an odd number, then $\dl(\Gamma_{n,k})=2$.
\end{proposition}

\proof
Let $k\ge 3$ be odd. Then $\Gamma_{n,k}$ is not bipartite and hence $\dl(\Gamma_{n,k})\geq 2$. To prove the reverse inequality we use the representation $\Gamma_{n,k} = P_{2n+1}(U)\sqcap C_k$, where $U=\{v_{2i+1}:\ 0\leq i \leq k\}$. Since $S_{P_{2n+1}}=\{v_1\}$ is a lmb, $S_{P_{2n+1}}\subseteq U$. Then by Theorem~\ref{thm:basis-under-U}, $\dl(\Gamma_{n,k}) = \dl(P_{2n+1}(U)\sqcap C_k)\leq \max\{\dl(P_{2n+1}), \dl(C_k)\}=2$. 
\qed

\subsection{Local metric basis in delivery services} 

Assume that a company wishes to assign codes to its customers such that the code of a customer uniquely determines its location. It is natural that the company is interested in making the length of the codes as short as possible. To design a graph theory model for this problem consider customers as the vertices of a(n) (edge-weighted) graph $G$. Vertices $u$ and $v$ are declared to be adjacent in $G$ if one of the following conditions is fulfilled:
\begin{itemize}
\item[(i)] there is no other customer on the $u,v$-geodesics;
\item[(ii)] the first letters of the family names of the customers $u$ and $v$ are the same.
\end{itemize}
To make the model more realistic, we also assign weights $w(uv)$ to the edges $uv$ of $G$ as follows. It the edge $uv$ is present solely because of (ii), then we set $w(uv) = \infty$. Otherwise (that is, if the edge $uv$ is present because of (i), or because of both (i) and (ii)), we set $w(uv)$ to be the real distance between the customers $u$ and $v$. 

Let $S$ be a lmb for $G$. Then, if $F$ is the first letter of the family name of a given customer $v$, then the company allocates the ordered pair $(F,r(v|S))$ to the customer $v$ as its code.
 
\begin{figure}[t!]
\centerline{\includegraphics[scale=.18]{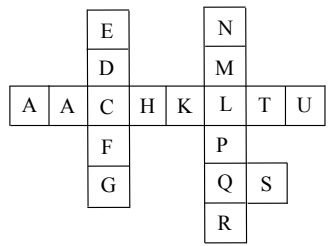}}
\caption{The robotic movement space (or the plan of customers' position).}
\label{fig:comp1}
\end{figure}

Let us compare the above model with the one suggested by Khuller et al.\ in their seminal paper on applications of the metric dimension~\cite{khuller-1996}. The model there was suggested to deal with robots' navigation in networks. As an example, consider the robotic movement space (or the plan of customers' position) as depicted in Fig.~\ref{fig:comp1}. Here all the distances in the weighted graph are set to $1$. 

\begin{figure}[t!]
\centerline{\includegraphics[scale=.23]{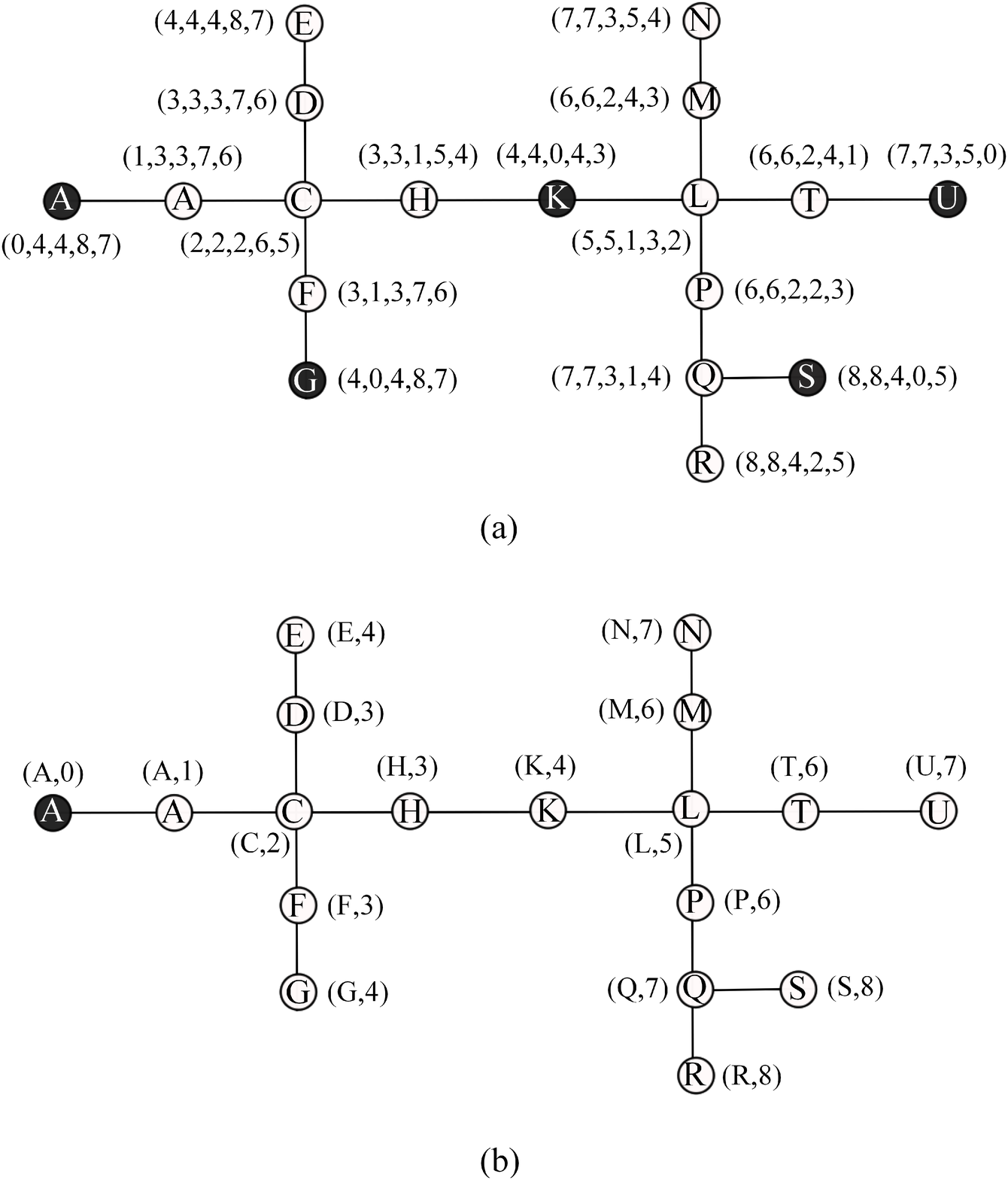}}
\caption{a) Earlier codes for the graph from Fig.~\ref{fig:comp1}. b) Present codes.}
\label{fig:comp4}
\end{figure}

In Fig.~\ref{fig:comp4}(a) the black vertices form a metric basis for the corresponding graph and the $5$-dimensional vectors next to vertices form the codes of the customers (locations) as proposed in the seminal model from~\cite{khuller-1996}. In Fig.~\ref{fig:comp4}(b) the new model is presented. Here $\{A\}$ is a lmb and the ordered pairs next to vertices are their codes. The second component of an ordered pair is formally a vector, but since its length is $1$, it is identified in the figure with the value of the component. 
 
Since $\dl(G)\leq {\rm dim}(G)$ holds for every connected graph $G$, cf.~\cite{okamoto-2010}, the proposed model is at least as compact as the earlier one. Furthermore, in many cases it is significantly shorter. For instance, $\dl(G)=1$ holds for every bipartite graph $G$, while ${\rm dim}(G)$ can be arbitrary large for such graphs.

\newpage

\newpage

\begin{thebibliography}{99}

\bibitem{anderson-2018}
  S.~E.~Anderson, S.~Nagpal, K.~Wash, 
  Domination in the hierarchical product and Vizing's conjecture,
  Discrete Math.\ 341 (2018) 20--24.

\bibitem{strong}
  G.~A.~Barrag\'an-Ram\'irez, J.~A.~Rodr\'iguez-Vel\'azquez,
  The local metric dimension of strong product graphs,
  Graphs Combin.\ 32 (2016) 1263--1278.

\bibitem{barriere-2009}
  L.~Barri\'ere, C.~Dafl\'o, M.~A.~Fiol, M.~Mitjana, 
  The generalized hierarchical product of graphs,
  Discrete Math.\ 309 (2009) 3871--3881.

\bibitem{chem1}
  G.~Chartrand, O.~R.~Oellermann, L.~Eroh, M.~A.~Johnson,
  Resolvability in graphs and the metric dimension of a graph,
  Discrete Appl.\ Math.\ 105 (2000) 99--113.

\bibitem{eliasi-2011}
  M.~Eliasi, A.~Iranmanesh, 
  The hyper-Wiener index of the generalized hierarchical product of graphs,
  Discrete Appl.\ Math.\ 159 (2011) 866--871.

\bibitem{fernau-2018}
  H.~Fernau, J.~A.~Rodr\'{i}guez-Vel\'{a}zquez, 
  On the (adjacency) metric dimension of corona and strong product graphs and their local variants: combinatorial and computational results,
  Discrete Appl.\ Math.\ 236 (2018) 183--202. 

\bibitem{gonzales-2018}
  A.~Gonz\'{a}lez, C.~Hernando, M.~Mora, 
  Metric-locating-dominating sets of graphs for constructing related subsets of vertices,
  Appl.\ Math.\ Comput.\ 332 (2018) 449--456. 
  
\bibitem{hakanen-2018}
  A.~Hakanen, T.~Laihonen, 
  On {$\{\ell\}$}-metric dimensions in graphs,
  Fund.\ Inform.\ 162 (2018) 143--160. 
  
\bibitem{Harary}
  F.~Harary, R.~A.~Melter, 
  On the metric dimension of a graph, 
  Ars Combin.\ (1976) 191--195.

\bibitem{hossein-2017}  
  S.~Hossein-Zadeh, A.~Iranmanesh, M.~A.~Hosseinzadeh, A.~Hamzeh, M.~Tavakoli, A.~R.~Ashrafi, 
  Topological efficiency under graph operations,
  J.\ Appl.\ Math.\ Comput.\ 54 (2017) 69--80.  

\bibitem{johnson-1993}
  M.~Johnson, 
  Structure-activity maps for visualizing the graph variables arising in drug design, 
  J.\ Biopharm.\ Stat.\ 3 (1993) 203--236.

\bibitem{kelenc-2018}
  A.~Kelenc, N.~Tratnik, I.~G.~Yero, 
  Uniquely identifying the edges of a graph: {T}he edge metric dimension,
  Discrete Appl.\ Math.\ 251 (2018) 204--220. 

\bibitem{khuller-1996}
  S.~Khuller, B.~Raghavachari, A.~Rosenfeld, 
  Landmarks in graphs, 
  Discrete Appl.\ Math.\ 70 (1996) 217--229.

\bibitem{chem2}
  A.~I.~Melker, S.~A.~Starovoitov, T.~V.~Vorobyeva,
  Classification of mini-fullerenes on graphs basis,
  Materials Phys.\ Mech.\ 20 (2014) 12--17.

\bibitem{okamoto-2010}
  F.~Okamoto, L.~Crosse, B.~Phinezy, P.~Zhang,
  The local metric dimension of a graph,
  Math.\ Bohem.\ 135 (2010) 239--255.

\bibitem{patt-2012}
  K.~Pattabiraman, P.~Paulraja, 
  Vertex and edge Padmakar-Ivan indices of the generalized hierarchical product of graphs,
  Discrete Appl.\ Math.\ 160 (2012) 1376--1384.

\bibitem{ecorona}
  Rinurwati, Slamin, H.~Suprajitno,
  General results of local metric dimensions of edge-corona of graphs, 
  Int.\ Math.\ Forum 11 (2016) 793--799.

\bibitem{rooted}
  J.~A.~Rodr\'{i}guez-Vel\'{a}zquez, C.~Garc\'{i}a G\'{o}mez, G.~A.~Barrag\'{a}n-Ram\'{i}rez, 
  Computing the local metric dimension of a graph from the local metric dimension of primary subgraphs,
  Int.\ J.\ Comput.\ Math.\ 92 (2015) 686--693.
 
\bibitem{Corona}
  J.~A.~Rodr\'{i}guez-Vel\'{a}zquez, G.~A.~Barrag\'{a}n-Ram\'{i}rez, C.~Garc\'{i}a G\'{o}mez, 
  On the local metric dimension of corona product graphs,
  Bull.\ Malays.\ Math.\ Sci.\ Soc.\ 39 (2016) S157--S173.  
 
\bibitem{saputro-2016}
  S.~W.~Saputro, 
  On local metric dimension of {$(n-3)$}-regular graph,
  J.\ Combin.\ Math.\ Combin.\ Comput.\ 98 (2016) 43--54.

\bibitem{Slater}
  P.~J.~Slater, 
  Leaves of trees, 
  Congr.\ Numer.\ 14 (1975) 549--559. 

\end{thebibliography}
\end{document}